\documentclass[]{amsart}
 \usepackage[]{amsmath, amsthm, amsfonts,amssymb}
 \usepackage[all]{xy}
 \usepackage[mathscr]{eucal}

 \newcommand {\C} {{\mathbb C}}
 \newcommand {\R} {{\mathbb R}}
 
 \newcommand {\Q} {{\mathbb Q}}
 \newcommand {\PP} {{\mathbb P}}

\newcommand {\OO} {{\mathcal O}}
\newcommand{\dt} {{\bullet}}

\newcommand {\spec} {{Spec\,}}
\newcommand {\G} {{\mathscr G}}
\newcommand {\ph} {\pi_1^{hodge}}
\newcommand {\IMHS} {Ind\text{-}MHS}
\newcommand {\T} {{\mathcal T}}
\newcommand {\E} {{\mathcal E}}
\newcommand {\V} {{\mathcal V}}

 \newtheorem{thm}[subsection]{Theorem}
 \newtheorem{cor}[subsection]{Corollary}
 \newtheorem{lemma}[subsection]{Lemma}
 \newtheorem{prop}[subsection]{Proposition}
 
 \newtheorem{rmk}[subsection]{Remark}
 \newtheorem{ex}[subsection]{Example}
  
 \newtheorem{conj}[subsection]{Conjecture}

\begin{document}

 \title{The Hodge theoretic fundamental group and its cohomology}
 
  \author{Donu Arapura}
\thanks{Partially supported by NSF}
 \address{Department of Mathematics\\
   Purdue University\\
   West Lafayette, IN 47907\\
   U.S.A.}
 \maketitle

From  the work of Morgan  \cite{morgan}, we know
that the fundamental group of
a complex algebraic variety carries a mixed Hodge structure,
which really means that a certain linearization
of it does.  This linearization, called the Mal\v{c}ev or pro-unipotent
completion destroys the group completely in some cases; for example if it were  perfect.
So a natural question is whether can one give a Hodge structure on a larger
chunk of the fundamental group. There have been a couple of approaches to this.
Work of Simpson \cite{simpson} continued by Kartzarkov, Pantev, Toen \cite{kpt1}
has shown that one has a weak  Hodge-like structure (essentially an action by $\C^*$ viewed
as a discrete group) on the entire pro-algebraic
completion of the fundamental group when the variety is smooth
projective. Hain \cite{hain3} refining his earlier work \cite{hain},
has shown that a Hodge structure of a more conventional sort
exists on the so called relative Mal\v{c}ev completions (under appropriate hypotheses).
In this paper, I want to propose a third alternative. I define a quotient
of the pro-algebraic completion called the Hodge theoretic
fundamental group $\ph(X,x)$ as the Tannaka dual to the category local systems
underlying admissible  variations of mixed Hodge structures on $X$, or in more
prosaic terms the inverse limit of Zariski 
closures of their associated monodromy representations.
This carries a nonabelian mixed Hodge structure in a sense that will 
be explained below.
The group $\ph(X,x)$ dominates the Mal\v{c}ev completion and the Hodge structures
are compatible, and I expect a similar statement for the relative completions.

The basic model here comes from arithmetic.  Suppose that $X$ is a variety
 over a field $k$  with  separable closure $\bar k$.
Let $\bar X = X\times_{\spec k}{\spec  {\bar k}}$. Suppose that $x\in X(k)$ is a
rational point, and that $\bar x\in X(\bar k)$ a geometric point lying over it.
Then there is an exact sequence of \'etale fundamental
 groups $$1\to \pi_1^{et}(\bar X,\bar x)\to \pi_1^{et}(X,\bar x)\to Gal(\bar
 k/k)\to 1$$
 The point $x$ gives a splitting, and so an action
 of $Gal(\bar k/k)$ on $\pi_1^{et}(\bar X,\bar x)$. This action will 
pass to the Galois cohomology of $\pi_1^{et}(\bar X,\bar x)$.
In the translation into Hodge theory, $\pi_1^{et}(\bar X,\bar x)$ is replaced
by $\ph(X,x)$ and
 the Galois group  by a certain universal Mumford-Tate group $MT$. The
action of $MT$ on $\ph(X,x)$ is precisely what I mean by a nonabelian mixed Hodge structure.
The cohomology $H^*(\ph(X,x),V)$ will carry induced mixed Hodge structures
for admissible variations $V$. In fact  there
is a canonical morphism $H^*(\ph(X,x),V)\to H^*(X,V)$. I call $X$ a Hodge theoretic $K(\pi,1)$
if this is  an isomorphism for all $V$.  Basic examples
of such spaces are abelian varieties, and  smooth affine
curves (modulo  \cite{hmpt}). I want to add that part of my motivation for this paper  is to test out some ideas which could be applied to motivic sheaves. So consequently certain constructions are phrased in more generality than is strictly
necessary for the present purposes.

My thanks to Roy Joshua for the invitation to the conference.
I would also like to thank Dick Hain,  Tony Pantev, Jon Pridham, and the referee for their
comments, and also Hain for informing me of
his recent work with Matsumoto, Pearlstein and Terasoma.

 \section{Review of Tannakian categories}

 In this section, I will summarize standard material from \cite{dm,deligne2, deligne-tan}.  Let
 $k$ be a field of characteristic zero.  Given an affine group scheme
 $G$ over $k$, we can express its coordinate ring as a directed union
 of finitely generated sub Hopf algebras $\OO(G)= \varinjlim A_i$.
 Thus we can, and will, identify $G$ with the pro-algebraic group
 $\varprojlim Spec\, A_i$ and conversely (see \cite[cor. 2.7]{dm} for
 justification).

 By a {\em tensor category} $\T$ over $k$,  I will mean a $k$-linear
 abelian category with bilinear tensor product making it into a
 symmetric monoidal (also called an ACU) category; we also require that
 the unit object ${\bf 1}$ satisfy $End({\bf 1})=k$. $\T$ is
 {\em rigid}  if it has duals.  The category $Vect_k$ of finite
 dimensional $k$-vector spaces is the key example of a rigid tensor
 category.  A {\em neutral Tannakian category} $\T$ over $k$ is a
 $k$-linear rigid tensor category, which possesses a faithful functor $F:\T\to
 Vect_k$, called a fibre functor, preserving all the structure.  The
 Tannaka dual $\Pi(\T,F)$ of such a category (with specified fibre
 functor $F$) is the group of tensor preserving automorphisms of $F$.
 In more concrete language, an element  $g\in\Pi(\T,F)$ consists of a collection
 $g_V\in GL(F(V))$, for each $V\in Ob\T$, satisfying the following
 compatibilities: 
 \begin{enumerate}
 \item[(C1)] $g_{V\otimes U} = g_V\otimes g_U$,
\item[(C2)]  $g_{V\oplus U} =
 g_V\oplus g_U$,
\item[(C3)]  the diagram
$$\xymatrix{
   F(V)\ar[r]^{g_V}\ar[d] & F(V)\ar[d] \\
   F(U)\ar[r]^{g_U} & F(U) } $$ commutes for every morphism $V\to U$.
 \end{enumerate}

 $\Pi=\Pi(\T,F)$ is (the group of $k$-points of) an affine group
 scheme. Suppose that
 $\T$ is generated, as a tensor category, by an object $V$ i.e. every
 object of $\T$  is finite sum of subquotients of tensor powers 
$$T^{m,n}(V) =  V^{\otimes m}\otimes V^{*\otimes n}.$$
Observe that for all $m,n\ge 0$
$$Hom_k(F({\bf 1}), F(T^{m,n}(V))) = T^{m,n}(F(V))$$

 \begin{lemma}\label{lemma:tannakaV}
Suppose that $\T$ is generated as a tensor category by $V$, then
\begin{enumerate}
\item $\Pi$  can be identified with the largest subgroup of $GL(F(V))$ that fixes all tensors in
the subspaces
$$Hom_\T({\bf 1}, T^{m,n}(V))\subset T^{m,n}(F(V))$$
\item $\Pi$  leaves invariant each subspace of $T^{m,n}(F(V))$ that corresponds to a subobject of $T^{m,n}(V)$.
\end{enumerate}
 In particular, $\Pi$ is an algebraic group.
 \end{lemma}

 \begin{proof}
   This is standard cf. \cite{dm}.
\end{proof}

  In general, 
$$\Pi(\T,F) = \varprojlim_{\T'\subset \T\text{ finitely generated}}
 Aut(F|_{\T'})$$ 
exhibits it as a pro-algebraic group and hence a group scheme.  

The following key example will explain our choice of notation.

 \begin{ex}
  Let $Loc(X)$ be the category of local systems (i.e. locally constant
  sheaves) of finite dimensional
$\Q$-vector spaces over a connected topological space $X$. This is a
neutral Tannakian category over $\Q$. For each $x\in X$, $F_x(L) =
L_x$ gives a fibre functor. The Tannaka dual  $\Pi(Loc(X),F_x)$ is isomorphic
to the rational pro-algebraic completion 
$$\pi_1(X,x)^{alg} =  \varprojlim_{\rho:\pi_1(X)\to GL_n(\Q)}
\overline{\rho(\pi_1(X,x))}$$
 of $\pi_1(X,x)$.
 \end{ex}

Given an affine
 group scheme $G$, let $Rep_\infty(G)$ (respectively $Rep(G)$) be the
 category of (respectively finite dimensional) $k$-vector spaces on
 which $G$ acts algebraically.  When $\T$ is  a neutral Tannakian 
category with a fibre functor $F$, the basic theorem of
 Tannaka-Grothendieck is that $\T$ is equivalent to $Rep(\Pi(\T,F))$. 
The role of a fibre functor  is similar to the role of
base points for the fundamental group of a connected
space. We can compare the groups at two base points by  choosing a
path between them. In the case of pair of fiber functors $F$ and $F'$,
a ``path'' is given by the   tensor 
isomorphism $p\in Isom(F,F')$ between
these functors. An element   $p\in Isom(F,F')$ determines an isomorphism
$\Pi(\T,F)\cong \Pi(\T,F')$ by $g\mapsto p g p^{-1}$. 
More canonically, one can define $\Pi(\T)$ as an ``affine group scheme in
$\T$'' independent of any choice of $F$ \cite[\S 6]{deligne2}. For our 
purposes, we can view  $\Pi(\T)$ as the Hopf algebra object which
maps to the coordinate ring $\OO(\Pi(\T,F))$ for each $F$.

Note that $\Pi$ is contravariant. That is,
given a faithful exact tensor functor $E:\T'\to \T$ between Tannakian
categories, we have an induced homomorphism $\Pi(\T,F)\to \Pi(\T',
F\circ E)$.

\section{Enriched local systems}

Before giving the definition of the Hodge theoretic fundamental group,
it is convenient to start with some generalities.
A theory of {\em enriched local systems} $E$  on the category smooth
complex varieties consist of

\begin{enumerate}
\item[(E1)]  an assignment of a neutral
$\Q$-linear Tannakian category $E(X)$ to every smooth variety $X$,
\item[(E2)]  a contravariant exact tensor pseudo-functor on the  category of
smooth varieties, i.e., a functor $f^*:E(Y)\to E(X)$ for each morphism
$f:X\to Y$ together with  natural isomorphisms for compositions,
\item[(E3)] faithful exact tensor functors $\phi:E(X)\to Loc(X)$ compatible with
  base change, i.e., a natural transformation of pseudo-functors $E\to
  Loc$,
\item[(E4)] a $\delta$-functor $h^\dt:E(X)\to E(pt)$ with natural isomorphisms
$\phi(h^i(L))\cong H^i(X,\phi(L))$. We also require there to be a
canonical map $p^*h^0(L)\to L$ corresponding the adjunction map on
local systems, where $p:X\to pt$ is the projection.
\end{enumerate}

By a weak  theory of enriched local systems, we mean something
satisfying (E1)-(E3). We have the following key examples.

\begin{ex}
  Choosing $E=Loc$ and $h^i=H^i$ gives the tautological example of a
   theory enriched of local systems.
\end{ex}

\begin{ex}
  Let $E(X)= MHS(X)$ be the category of admissible variations of mixed
  Hodge structure \cite{kashiwara,sz} on $X$. This carries a forgetful
  functor $MHS(X)\to Loc(X)$. Let $h^i(L)=H^i(X,L)$ equipped with the mixed
  Hodge structure constructed by Saito \cite{saito1, saito}. Here
  $p^*h^0(L)$ is the invariant part $L^{\pi_1(X)}\subset L$. This
  can be seen to give a sub variation of MHS by restricting to
  embedded curves and applying \cite[4.19]{sz}. Thus $MHS(X)$ is
  a theory of enriched local systems.
\end{ex}

Many more examples of  theories of enriched local systems can be obtained by taking subcategories of the ones
above.

\begin{ex}
 The category $E(X)= HS(X)\subset MHS(X)$ of direct sums of pure  variations of Hodge
 structure of possibly different weights.
\end{ex}

\begin{ex}
  The category $E(X) = UMHS(X)\subset MHS(X)$ of unipotent variations of mixed Hodge structure.
\end{ex}

\begin{ex}
  Finally the category of tame motivic local systems constructed in
  \cite[sect 5]{arapura} is, for the present, only a weak  theory of enriched
  local systems. However, the category defined by systems of
  realizations [loc. cit.] can be seen to be an enriched theory in the
  full sense.
\end{ex}

Given a theory of enriched  local systems $(E,\phi, h^\dt)$, let
$\phi(E(X))$ denote the Tannakian subcategory of $Loc(X)$ generated by the image
of $E(X)$. So $\phi(E(X))$ is the full subcategory whose
objects are sums of subquotients of objects in the image of $\phi$.
Set $\pi_1^E(X,x) = \Pi(\phi(E(X)), F_x)$, where $F_x$ is the fibre functor
associated to a base point $x$. More explicitly, this is the inverse
limit of the Zariski closures of monodromy representations of objects of
$E(X)$. It follows that the isomorphism class of $\pi_1^E(X,x)$ as a
group scheme is
independent of $x$. However, certain additional structure will
depend on it.

Let $\kappa:E(pt)\to E(X)$ and  $\psi: E(X)\to E(pt)$ be given by
$p^*$ and $i^*$ respectively, where  $p:X\to pt$ and $i:pt\to X$ 
are the projection and inclusion of $x$. We also have $\phi:E(X)\to
\phi(E(X))$. These functors yield a diagram
$$\pi_1^E(X,x)\to \Pi(E(X),x)\rightleftarrows \Pi(E(pt))$$
where $\Pi(E(X),x) =\Pi(E(X), F_x)$ and $\Pi(E(pt)) =
\Pi(E(pt),pt)$. The diagram is clearly canonical in the sense
that a morphism $f:(X,x)\to (Y,y)$ of pointed varieties gives rise to
a larger commutative diagram
$$
\xymatrix{
 \pi_1^E(X,x)\ar[r]\ar[d] & \Pi(E(X), x)\ar@<1ex>[r]\ar[d] & \Pi(E(pt))\ar[l]\ar[d]^{=} \\ 
 \pi_1^E(Y,y)\ar[r] & \Pi(E(Y), y)\ar@<1ex>[r] & \Pi(E(pt))\ar[l]
}
$$

\begin{thm}\label{thm:PiMT}
  The sequence $$1\to \pi_1^E(X,x)\to \Pi(E(X),x)\rightleftarrows \Pi(E(pt))\to 1$$
  is split exact. Therefore there is a canonical isomorphism
$ \Pi(E(X),x)\cong \Pi(E(pt))\ltimes \pi_1^E(X,x)$. 
\end{thm}

\begin{proof}
  Since $\psi\circ\kappa = id$, the induced homomorphisms $\Pi(E(pt))\to
  \Pi(E(X),x)\to \Pi(E(pt))$ compose to the identity.  The injectivity of
  $\pi^E_1(X,x)\to \Pi(E(X),x)$ follows from \cite[2.21]{dm}.

  Therefore it remains to check exactness in the middle.  An element
  of $ im[\pi_1^E(X,x)\to \Pi(E(X),x)]$ is given by a collection of
   elements $g_V\in GL(\phi(V)_x)$, with
  $V\in ObE(X)$, satisfying (C1)-(C3) such that
 \begin{equation}\label{eq:fUV}
    g_U = \alpha_x\circ g_V\circ\alpha_x^{-1}  
  \end{equation} 
  holds for  any isomorphism $\alpha:\phi(V)\cong \phi(U)$.
 While  an element of $ker[\Pi(E(X),x)\to \Pi(E(pt))]$ is given by a
  collection $\{g_V\}$ such that
  $g_V=I$ for any object in the image of $\kappa$. 
If $V$ is in the image of $\kappa$, the underlying local system is
trivial. This implies that $g_V = \alpha_x\circ\alpha_x^{-1}=I$, if
(\ref{eq:fUV}) holds.
Thus, we
  have $$im[\pi_1^E(X,x)\to \Pi(E(X),x)]\subseteq
  ker[\Pi(E(X),x)\to \Pi(E(pt)]$$

  Conversely, suppose that $\{g_{V}\}$ is an element of the kernel on the
  right.  Given an isomorphism $\alpha:\phi(V)\to
  \phi(U)$,  we have to verify (\ref{eq:fUV}). Let
  $H=\kappa(h^0(V^*\otimes U))$, then $g_H=1$ by assumption. 
 $H$ gives a subobject of $V^*\otimes
  U$ which maps to the invariant part  $\phi(V^*\otimes U)^{\pi_1(X)}$
   of the local system $\phi(V^*\otimes U)$. This follows from our
   axiom (E4). Therefore, $\alpha$ gives
   a section of $\phi(H)$. The evaluation morphism
   $ev:(V^*\otimes U)\otimes V\to U$ restricts to give a morphism
   $H\otimes V\to U$. We claim that the  diagram
$$
\xymatrix{
 \phi(V)_x\ar[r]_<<<<{\alpha\otimes I}\ar[d]^{g_V}\ar@/^1pc/[rr]^{\alpha} & \phi(H)_x\otimes \phi(V)_x\ar[r]_<<<<{ev}\ar[d]^{I\otimes g_V} & \phi(U)_x\ar[d]^{g_U} \\ 
 \phi(V)_x\ar[r]^<<<<{\alpha\otimes I}\ar@/_1pc/[rr]_{\alpha} & \phi(H)_x\otimes \phi(V)_x\ar[r]^<<<<{ev} & \phi(U)_x
}
$$

commutes.
The commutativity of the square on the left is clear.
For the commutativity on right, apply (C1) and (C3) and the fact that
 $g_H = I$.
Equation (\ref{eq:fUV}) is now proven.
\end{proof}

Since any representation of an affine group scheme is locally finite, it follows that
$Rep_\infty(E(pt))$ can be identified with the category of ind-objects
$Ind\text{-}E(pt)$.
We will often refer to an object of this category as an $E$-structure.
To simplify notation, we generally will not distinguish between $V$ and $\phi(V)$.
By a {\em nonabelian  $E$-structure}, we will mean  an affine group scheme $G$
over $\Q$ with an algebraic action of $\Pi(E(pt))$. Equivalently,
$\OO(G)$ possesses an $E$-structure compatible with the Hopf algebra operations.
 A morphism of  nonabelian  $E$-structures is a
homomorphism of group schemes commuting with the $\Pi(E(pt))$-actions.
The previous theorem yields a nonabelian $E$-structure on  $\pi_1^E(X,x)$ which
is functorial in the category of pointed varieties. When $X$ is connected, 
for any two base points $x_1,x_2$,
 $\pi_1^E(X,x_1)$ and  $\pi_1^E(X,x_2)$ are isomorphic as group schemes although the
$E$-structures need not be the same.

\begin{lemma}\label{lemma:NES}
  An algebraic group $G$ carries a nonabelian
  $E$-structure if and only if there exists a finite dimensional
  $E$-structure $V$ and an embedding $G\subseteq GL(V)$, such that $G$ is
  normalized by the action of $\Pi(E(pt))$.
  A  nonabelian $E$-structure is an inverse limit of algebraic groups
  with $E$-structures. 
\end{lemma}

\begin{proof}
  
  Suppose that $G$ is an algebraic group with an $E$-structure.
  Then $\Pi= \Pi(E(pt))$ acts on $G$; denote this left action by
  ${}^mg$. We also have a left action of $\Pi$ on $\OO(G)$, written in
  the usual way, such that
  \begin{equation}
    \label{eq:mgf}
  m(gf)= ({}^mg)( mf)    
  \end{equation}
  for $m\in \Pi$, $g\in G$ and $f\in
  \OO(G)$. This gives an action of $\Pi\ltimes G$ on $\OO(G)$.
  Let $V'\subset \OO(G)$ be a subspace spanned by a finite set of algebra
  generators.  Let $V\subset \OO(G)$
  be  the smallest $\Pi\ltimes G$-submodule containing $V'$.  By standard
  arguments, $V$ is  finite dimensional and faithful as a $G$-module. So we get $G\subseteq
  GL(V)$. Equation (\ref{eq:mgf}) implies that $G$ is normalized by
  $\Pi$.  This proves one direction. The converse is clear.

For the second statement write $\OO(G)$ as direct limit $\OO(G)
  =\varinjlim V_i$ of finite dimensional $E$-structures. Let
  $A_i\subset \OO(G)$ be the smallest Hopf
  subalgebra containing $V_i$. This is an $E$-substructure. 
  So we have $G= \varprojlim A_i$ which gives the desired conclusion.
\end{proof}

An {\em $E$-representation} of a  nonabelian $E$-structure $G$
is a representation on an $E$-structure $V$ such that  (\ref{eq:mgf})
holds. The adjoint representation gives
a canonical example of an $E$-representation.
The above lemma says
that every finite dimensional nonabelian $E$-structure has a
faithful $E$-representation.  The following is
straightforward and was used implicitly already.

\begin{lemma}\label{lemma:semidir}
  There is an equivalence between the tensor category of
  $E$-representations of a nonabelian $E$-structure $G$ and the
  category of representations of the semidirect product $\Pi(E(pt)) \ltimes
  G$.
\end{lemma}

\begin{cor}
  In the notation of theorem~\ref{thm:PiMT}, $E(X)$ is equivalent to
  the category of finite dimensional $E$-representations of $\pi_1^E(X,x)$.
\end{cor}

\section{Nonabelian Hodge structures}

 The category $MHS=MHS(pt)$ (respectively $HS=HS(pt)$) is just the category of 
 rational graded polarizable rational mixed (respectively pure) Hodge
 structures \cite{deligne1}; its Tannaka dual
will be called the universal (pure) Mumford-Tate group and will be denoted by
$MT$ ($PMT$).  The Tannaka dual of the category of real Hodge structures $\R\text{-}HS$ is just
Deligne's torus ${\mathbb S} = Res_{\C/\R}\C^*$. So the obvious functor $\R\text{-}HS\to HS\otimes \R$
yields an embedding ${\mathbb S}(\R) \hookrightarrow PMT(\R)$.  In more concrete terms, $ (z_1,z_2)\in
 {\mathbb S}(\C) = \C^*\times \C^*$ acts by multiplication by
 $z_1^p\bar z_2^{q}$ on the $(p,q)$ part of a pure Hodge structure.
Since $HS$ is semisimple, $PMT$
is pro-reductive. 
The inclusion $HS\subset MHS$ gives a homomorphism
$MT\to PMT$. We have a section $PMT\to MT$ induced by the functor
$V\mapsto Gr^W(V) = \oplus Gr_i^W(V)$.  Thus $MT$ is a semidirect
product of $PMT$ with $\ker[MT\to PMT]$. The kernel is pro-unipotent since it
acts trivially on $W_\dt V$ for any $V\in MHS$.

$Rep_\infty(MT)$ is the category $Ind\text{-}MHS$ of direct limits of mixed
Hodge structures. Given an object $V= \varinjlim V_i$ in this category,
we can extend the Hodge and weight filtrations by $F^p V =
\varinjlim F^p V_i$ and $W_k V = \varinjlim W_k V_i$.

 Set $\ph= \pi_1^E$ for $E=MHS$. 
So this is the inverse limit of Zariski closures of monodromy
representations of variations of mixed Hodge structures.
The key definition is:
\begin{enumerate}

\item[(NH1)] 
A  nonabelian mixed (respectively pure) Hodge structure, or simply an $NMHS$ (or an
$NHS$), is an affine group scheme $G$
over $\Q$ with an algebraic action of $MT$ (respectively $PMT$). 
\end{enumerate}

A morphism of these objects is a
homomorphism of group schemes commuting with the $MT$-actions. 
A Hodge representation is an $MT$-equivariant representation; it is 
the same thing as an $E$-representation for $E=MHS$.
The coordinate ring of an $NMHS$ is a Hopf algebra in $\IMHS$.
Let us recapitulate the results of the previous section in the present
setting.
\begin{itemize}
\item  $\ph(X,x)$ has an $NMHS$ which is
functorial in the category of smooth pointed varieties (and this structure
will usually depend on the choice of $x$).
\item An algebraic group $G$ admits an $NMHS$ if and only if it has a
  faithful representation to  the general linear group of a mixed 
  Hodge structure for which $MT$ normalizes $G$.
\item Admissible variations of MHS correspond to Hodge representations
  of $\ph(X,x)$. 
\end{itemize}

The notion of a nonabelian mixed Hodge structure is  fairly weak,
although sufficient for some of the main results of this paper.
At this point it is not clear  what the optimal set
of axioms should be.
We would like to spell out some further conditions which will hold in our basic example
$\pi_1^{hodge}$.

\begin{enumerate}

\item[(NH2)] 
An NMHS $G$ satisfies (NH2) or  has nonpositive weights if $W_{-1}\OO(G)=0$.
\end{enumerate}

\begin{rmk}
To see why this is``nonpositive'', observe that  if $V$ is an MHS with $W_1V=0$, then
$W_{-1}\OO(V) = W_{-1}Sym^*(V^*)=0$. (My thanks to the referee for pointing out that the weights gets flipped.)
\end{rmk}

The significance of this condition  is explained by the following:

\begin{lemma}\label{lemma:NH2}
An NMHS $G$ has nonpositive weights if and only if the left and right actions of $G$ on $\OO(G)$ preserves  the weight filtration.
\end{lemma}

\begin{proof}
Note that the left and right $G$-actions preserves the weight filtration if and
only if
\begin{equation}
  \label{eq:muW1}
\mu^*(W_i\OO(G))\subseteq\OO(G)\otimes W_i\OO(G)  
\end{equation}
\begin{equation}
  \label{eq:muW}
\mu^*(W_i\OO(G))\subseteq W_i\OO(G)\otimes\OO(G)  
\end{equation}

Comultiplication $\mu^*:\OO(G)\to \OO(G)\otimes \OO(G)$ is a morphism of $Ind$-MHS.
If  $G$ has nonpositive weights, then
$$\mu^*(W_i\OO(G))\subseteq   \sum_{j\ge 0}  W_j\OO(G)\otimes W_{i-j}\OO(G)\subseteq \OO(G)\otimes W_i\OO(G)$$
 \eqref{eq:muW} follows by symmetry.

Suppose $W_{-1}\OO(G)\not=0$ and that \eqref{eq:muW1}, and \eqref{eq:muW} hold.
Let $f\in W_{-1}\OO(G)$ be a nonzero element, and let $n>0$ be the largest
integer such that $f\in W_{-n}$. Suppose that
$\mu^*(f) = \sum g_i\otimes h_i$. By \eqref{eq:muW},
$$\sum g_i\otimes h_i\equiv 0 \mod \OO(G)/W_{-n}\otimes \OO(G)$$
Therefore all $g_i\in W_{-n}$. By a similar argument $h_i\in W_{-n}$.
Therefore $\mu^*(f)\in W_{-2n}(\OO(G)\otimes \OO(G))$. Since morphisms
of MHS (and therefore Ind-MHS) strictly preserve weight filtrations \cite{deligne1},
$\mu^*(f) =\mu^*(f')$ for some $f'\in W_{-2n}$. But $\mu^*$ is injective because $\mu$ is dominant. Therefore $f=f'\in W_{-2n}$ which is a contradiction.
\end{proof}

\begin{lemma}
 If $G$ satisfies (NH2), then there is a unique maximal pure quotient $G^{pure}$.
\end{lemma}

\begin{proof}
Since $G$ satisfies (NH2) it acts on the pure Ind-MHS 
$Gr^W\,\OO(G)$ on the left.  As  a group we  take $G^{pure}$ to be
the image of $G$ in $Aut(Gr^W\,\OO(G))$.
To get the finer structure, we apply lemma~\ref{lemma:NES}, to
write $G=\varprojlim G_i$ with $G_i\subset GL(V_i)$ where $V_i\subset
\OO(G)$ are mixed Hodge structures such that $G_i$ is normalized by $MT$.
 $G$ preserves $W_\dt V_i$ by assumption. Then 
$$G^{pure}= \varprojlim im[G_i\to GL(Gr^W\,V_i)]$$
This group carries a nonabelian pure Hodge structure
since each group of the limit does.
By construction, there is a surjective morphism $G\to G^{pure}$.

Suppose that $G\to H$ is another pure quotient. Then $G$ will act on
$\OO(H)$ through this map. The image of $G$ in $Aut(\OO(H))$ is precisely $H$.
We have a $G$-equivariant morphism of Ind-MHS
$\OO(H)\subset \OO(G)$. By purity $\OO(H) = Gr^W(\OO(H))\subset Gr^W(\OO(G))$,
 which shows that the $G$-action on $\OO(H)$ factors through $G^{pure}$.
Therefore the homomorphism $G\to H$ also factors through this.
\end{proof}

Before explaining  the next condition, we recall that some
standard definitions. Fix a real algebraic group $G$. We have a
conjugation $g\mapsto \bar g $ on the group of complex points $G(\C)$,
whose fixed points are exactly $G(\R)$.  A {\em Cartan involution} $C$ of
$G$ is an algebraic involution of $G(\C)$ defined over $\R$, such that
 the group of fixed points of $g\mapsto C(\bar g)=\overline{C(g)}$ is compact
for the classical topology and has a point in every component. 
An involution of a pro-algebraic group is Cartan if it descends to a Cartan
involution in the usual sense on a cofinal system of finite
dimensional quotients.
 Recall that Deligne's torus ${\mathbb S} = Res_{\C/\R}\C^*$ embeds
into $PMT$ in such a way that  $ (z_1,z_2)\in
 {\mathbb S}(\C) = \C^*\times \C^*$ acts by multiplication by
 $z_1^p\bar z_2^{q}$ on the $(p,q)$ part of a pure Hodge structure.
 It will be convenient to say that a group is reductive (or pro-reductive) when its
connected component of the identity is.

\begin{enumerate}

\item[(NH3)] A NMHS $G$ satisfies (NH3) or is S-polarizable if it has nonpositive weights, and
 the action $C$ of $(i,i)\in {\mathbb  S}(\C)$ on $G^{pure}$ gives a 
Cartan involution.    
\end{enumerate}
The ``S'' stands for Simpson, since
this condition  is related to (and very much inspired by) the notion of
a pure  nonabelian Hodge structure 
introduced by him \cite[p. 61]{simpson}.
 Simpson has shown that the
pro-reductive completion of the fundamental group of a smooth
projective variety carries a pure nonabelian Hodge structure in his
sense. These ideas have been developed further in  \cite{kpt1, kpt2}.
Although there is no direct relation between these notions of
nonabelian Hodge structure, there
are a number of close parallels (e.g. lemma
\ref{lemma:hodgetype} below holds for both). 
The meaning of (NH3) is explained by the following:

\begin{lemma}\label{lemma:NPHS}
  Let $G$ be an algebraic group with an NMHS. Then $G$ is $S$-polarizable if and only
  if there  exists  a pure Hodge structure $V$ and an embedding
  $G^{pure}\subseteq GL(V)$, such that $G^{pure}$ is normalized by the action of
  $PMT$ and such that $G^{pure}$
  preserves a polarization on $V$. Under these conditions, $G^{pure}$ is reductive.
\end{lemma}

\begin{proof}
After replacing $G$ by $G^{pure}$, we may assume  $G=G^{pure}$ is already pure.
  Given an embedding   $G\subseteq GL(V)$ as above,
  choose a $G$-invariant polarization $(\,,\,)$ on
  $V$.  The image $W$ of $(i,i)\in {\mathbb S}(\C)$ in $GL(V)$ is
  nothing but the Weil operator for the Hodge structure on $V$.
  Therefore $\langle u,v\rangle = (u,W\bar v)$ is positive definite
  Hermitian.  The group of fixed points under $\sigma(g)=W^{-1}\bar
  gW$ is easily seen to preserve $\langle \, ,\,\rangle$, so it is
  compact.  In other words, $Cg=W^{-1} gW$ is a Cartan involution.
  Conversely, if $C$ is a Cartan involution then a $G$-invariant
  polarization can be constructed by averaging an existing
  polarization over the Zariski dense compact group of $\sigma$-fixed
  points.  When these conditions are satisfied the reductivity of $G$
  follows from the existence of a Zariski dense compact subgroup.
\end{proof}

\begin{cor}
  The underlying group of a pure S-polarizable nonabelian Hodge structure
   is  pro-reductive.
\end{cor}

\begin{cor}
If $G$ is S-polarizable, then  $G^{pure}=G^{red}$ is the maximal
pro-reductive quotient of $G$. In particular, if $G$ is pro-reductive then
$G= G^{pure}$ is pure.
\end{cor}

\begin{proof}
 We have an exact sequence of pro-algebraic groups 
$$1\to U\to G\to G^{pure}\to 1$$
where $U$ is simply taken to be the kernel.  This can
also be described as above by $$U = \varprojlim \ker[G_i\to
GL(Gr^W\,V_i)]$$
We can see from this that $U$ is pro-unipotent.  By the previous
corollary, $G^{pure}$ is pro-reductive. Thus $U$ is the
pro-unipotent radical and $G^{pure}=G^{red}$ is the maximal
pro-reductive quotient. This forces  $G= G^{pure}$ if $G$ is pro-reductive.
\end{proof}

Recall that Simpson \cite[p 46]{simpson} defined a real algebraic
group $G$ to be of {\em Hodge type} if $\C^*$ acts on $G(\C)$ such that
$U(1)$ preserves the real form and $-1\in U(1)$ acts as a Cartan
involution.  Such groups are reductive, and also subject to a number
of other restrictions [loc. cit.]. For example, $SL_n(\R)$ is not of Hodge type
when $n\ge 3$.

\begin{lemma}\label{lemma:hodgetype}
  If an algebraic group admits a pure S-polarizable nonabelian Hodge
  structure then it is of Hodge type.
\end{lemma}

\begin{proof}
  Choose an embedding $G\subset GL(V)$ as in lemma~\ref{lemma:NPHS}.
  The Hodge structure on $V$ determines a representation ${\mathbb
    S}(\C)\to GL(V_\C)$. The group $G(\C)$ is stable under conjugation
  by elements of $PMT(\C)\supset \mathbb{S}(\C)$. Embed
  $\C^*\subset {\mathbb S}(\C)$ by the diagonal.  Then $-1\in \C^*$
  acts trivially on $G(\C)$. Therefore the $\C^*$ action factors through
  $\C^*/\{\pm1\} \cong \C^*$. To see that $U(1)/\{\pm 1\}$ preserves
  $G(\R)\subset GL(V_\R)$, it suffices to note that the image of
  $e^{i\theta}\in U(1)$ in $GL(V_\C)$, which acts on $V^{pq}$ by
  multiplication by $e^{(p-q)i\theta}$, is a real operator.  Under the
  isomorphism $U(1)\cong U(1)/\{\pm 1\}$, $-1$ on the left corresponds to the
  image of $i$ on the right.  Thus $-1$ acts by a Cartan involution on $G(\C)$.
 
\end{proof}

\begin{thm}\label{thm:ph}
 Given a smooth variety $X$,
  $\ph(X,x)$ carries an S-polarizable nonabelian Hodge mixed structure. 
The category of Hodge representations of
  $\ph(X,x)^{red}$ is equivalent to the category of pure variations
  of Hodge structure on $X$.
\end{thm}

\begin{proof} 
   For $V\in MHS(X)$, let $\ph(\langle V\rangle,x)$ denote the Zariski
  closure of the monodromy representation of $\pi_1(X,x)\to GL(V_x)$.
   Let $MT(V)\subset GL(V_x)$ denote
  the Tannaka dual of the sub tensor category of $MHS$ generated by
  $V_x$. A slight modification of theorem~\ref{thm:PiMT} together with
  lemma~\ref{lemma:NES} shows that $\ph(\langle V\rangle,x)$ is normalized by
  $MT(V)$. We can also see this  directly. The group $\ph(\langle V\rangle,x)$ is  characterized as the group of
  automorphisms that fix all monodromy invariant tensors
  $T^{m,n}(V_x)^{\pi_1(X,x)}$. While $MT(V)$ 
  leaves all sub MHS of  $T^{m,n}(V_x)$ invariant by lemma~\ref{lemma:tannakaV}.  Let
  $g\in MT(V)$ and let $\gamma\in \ph(\langle V\rangle,x)$, then it is enough to see
  that $g^{-1}\gamma g$ fixes every  tensor in
  $T^{m,n}(V_x)^{\pi_1(X,x)}$. This space is a sub MHS of
  $T^{m,n}(V_x)$ by \cite[4.19]{sz}. Therefore $T^{m,n}(V_x)^{\pi_1(X,x)}$ is invariant
  under $g$ (although it need not fix elements pointwise). This shows
  that $g^{-1}\gamma g$ fixes the elements of this space as
  claimed. Therefore $\ph(\langle V\rangle,x)$ carries a NMHS.

  Since the weight filtration of $V$ is a filtration by local systems,
  $\ph(\langle V\rangle,x)$ preserves $W_\dt V_x$. So it satisfies (NH2) by lemma \ref{lemma:NH2}. Let $\tilde \pi_1^{hodge}(V)$ be
  the image of $\ph(\langle V\rangle,x)$ in $GL(Gr^W \, V_x)$. This can be
  identified with the Zariski closure of the monodromy representation
  of the pure variation of Hodge structure $Gr^W V$. This pure
  variation is polarizable by definition of
  admissibility \cite{sz, kashiwara}. Therefore $Gr^W \, V_x$ possesses a $\tilde
  \pi_1^{hodge}(V)$ invariant polarization. Consequently
  $\ph(\langle V\rangle,x)$ satisfies (NH3) by lemma~\ref{lemma:NPHS}. Moreover $\tilde \pi_1^{hodge}(V)=
  \ph(\langle V\rangle,x)^{red}$ is the Tannaka dual to the subcategory of
  $Loc(X)$ generated by the local system $Gr^W \, V$.  Putting this
  all together, we see that $$\ph(X,x)=\varprojlim_V \,
  \ph(\langle V\rangle,x)$$
  satisfies (NH3), and that the
  Tannaka dual to $HS(X)$ is $PMT\ltimes \ph(X,x)^{red}$.
  Lemma~\ref{lemma:semidir} implies that $HS(X)$ is equivalent
  to the category of Hodge representations of $\ph(X,x)^{red}$.
\end{proof}

\section{Nonabelian variations}

The goal of this section is to give a characterization of $\ph(X,x)$. To this end,
we introduce  the category of {\em nonabelian variation of mixed Hodge structures}  over $X$,
by which we mean the opposite of the category of Hopf algebras in Ind-$MHS(X)$. 
Given such a Hopf algebra $A$, we denote the corresponding nonabelian variation by the {\em symbol}
$Spec\, A$.  For any $x\in X$, $Spec\, A_x $ can be understood in the usual sense, and this
is an NMHS.  Basic examples are given as follows.

\begin{ex}\label{ex:VMHS1}
  Any object $V\in MHS(X)$ can be identified with the nonabelian variation
$Spec\, (Sym^*(V^*))$.
\end{ex}

By applying the forgetful functor $MHS(X)\to Loc(X)$, we can see that any nonabelian
variation $Spec\, A$ carries a monodromy action of $\pi_1(X,x)\to Aut(G_x)$,
where $G_x =Spec\, A_x$. 
A nonabelian variation of mixed Hodge structures will be called {\em inner}
if  the monodromy action lifts to a homomorphism $\pi_1(X,x)\to G_x$.
The examples of \ref{ex:VMHS1} are rarely inner. However,  an ample
supply of such examples is given by the following.

\begin{ex}\label{ex:VMHS2}
  For $V\in MHS(X)$,  we have a nonabelian variation $\ph(\langle V\rangle)$
whose fibres are $\ph(\langle V\rangle, x)$ by \cite[\S 6]{deligne2}. To construct this directly, note that
we can realize the coordinate ring of  $\ph(\langle V\rangle, x)$ as a quotient of
$\OO(GL(V_x))$ by a Hopf ideal $\sum f_k\OO(GL(V_x))$. The generators $f_k$
can be regarded as sections of  
$$\mathcal{R}= \bigoplus_{i,j\ge 0} Sym^i(V^*)\otimes \det(V^*)^{-j} $$
Thus we can define  $\ph(\langle V\rangle)$ as $Spec$ of the Ind-$MHS(X)$ Hopf algebra
$\mathcal{R}/(\sum_k f_k \mathcal{R})$.
This is inner
since the monodromy is given by homomorphism $\pi_1(X,x)\to \ph(\langle V\rangle, x)$.
\end{ex}

Let $\ph(X)$ be the inverse limit of  $\ph(\langle V\rangle)$ over $V\in MHS(X)$. The fibres
are $\ph(X,x)$. This is the universal inner nonabelian variation in the following sense:

\begin{prop}\label{prop:NVMHS}
If $G$ is  an inner nonabelian variation of mixed Hodge structure  over $X$,
with monodromy given by a homomorphism
$\rho:\pi_1(X,x)\to G_x$. Then $\rho$ extends to a morphism 
$\ph(X)\to G$ of nonabelian variations. 
\end{prop}

\begin{proof}
Let $G=Spec\, A$. Since $A\in$ Ind-$MHS(X)$, $\ph(X,x)$ will act on it.
By an  argument similar to the proof of lemma~\ref{lemma:NES},
we can write $A$ as a direct limit of finitely generated Hopf algebras
with $\ph(X,x)$-action. So that $G_x$ becomes an inverse limit of algebraic groups 
carrying inner nonabelian variations. 
Thus we can assume that $G_x$ is an algebraic group.
   By standard techniques
we can find a finite dimensional faithful (left) $G$-submodule $V\subset \OO(G)$.
After replacing this with the span of the $\ph(X,x)$-orbit, we can assume that
$V$ is stable under $\ph(X,x)$. Therefore $V$ corresponds to a variation of
mixed Hodge structure. By assumption, the image of $\pi_1(X,x)$ in $GL(V_x)$ lies
in $G$. This implies that $G$ contains $\ph(\langle V\rangle,x)$ and that $\rho$
factors through it.
\end{proof}

\section{Unipotent  and Relative Completion}

Morgan \cite{morgan} and Hain \cite{hain} have shown that the
pro-unipotent completion of the fundamental group of a smooth variety
carries a mixed Hodge structure. We want to compare this with our
nonabelian Hodge structure. We start by recalling  some standard facts
from group theory (c.f. \cite{hz}, \cite[appendix A]{quillen}).
Fix a finitely generated group $\pi$. Then
\begin{enumerate}
\item[(a)]  $\Q[\pi]$, and its quotients by powers of the augmentation
  ideal $J$, carry  Hopf algebra structures with comultiplication
$\Delta(g)= g\otimes g \mod J^r$.
\item[(b)] A finite dimensional $\Q[\pi]$-module is unipotent if and only
  if it factors through some power of $J$. (The smallest power will be called
  the index of unipotency).
\item[(c)] The set of group-like elements $$G_r(\pi)= \{f\in
  \Q[\pi]/J^{r+1}\mid \Delta(f) = f\otimes f,\> f\equiv 1\mod J\}$$
  forms a group under   multiplication. This is a unipotent algebraic group.
\item[(d)] The Lie algebra of $G_r(\pi)$ can be identified with the Lie
  algebra of primitive elements $$
  {\mathscr G}_r(\pi) = \{f\in
  \Q[\pi]/J^{r+1}\mid \Delta(f) = f\otimes 1+1\otimes f\}$$
with bracket given by commutator.

\item[(e)] The exponential map gives a bijection of sets
${\mathscr G }_r(\pi)\cong G_r(\pi)$. The Lie algebra and group structures determine
each other via the Baker-Campbell-Hausdorff formula.
\item[(f)] 
  $\Q[\pi]/J^{r+1}$ is isomorphic to a quotient of the universal enveloping  algebra
   of $\mathscr{G}_r(\pi)$ by  a power of its augmentation ideal.
\end{enumerate}

Let $ULoc(X)$ ($U_rLoc(X)$) denote the category of local
systems with unipotent monodromy (with index of unipotency at most
$r$). The category $U_rLoc(X)$ can be identified with the
category of $\Q[\pi_1(X,x)]/J^{r+1}$-modules. 
We note that  this category has a tensor product:
$\Q[\pi_1(X,x)]/J^{r+1}$ acts on the usual tensor product of
representations $U\otimes_\Q V$ through $\Delta$.
With this structure, the category of $\Q[\pi_1(X,x)]/J^{r+1}$-modules
is Tannakian. Its   Tannaka dual
$\pi_1(X,x)^{un_{r}}$ is
isomorphic to the group $G_r(\pi_1(X,x))$ above, and the Tannaka dual
$\pi_1(X,x)^{un}$ of $ULoc(X)$, is the inverse limit of these groups.

We need to impose Hodge structures on these objects. 

\begin{lemma}\label{lemma:nmhspi1}
  There is a bijection between
  \begin{enumerate}
  \item The set of nonabelian mixed Hodge structures on $G_r(\pi)$.
  \item The set of mixed Hodge structures on $\mathscr{G}_r(\pi)$
    compatible with Lie bracket.
  \item The set of mixed Hodge structures on $\Q[\pi]/J^{r+1}$
    compatible with the Hopf algebra structure.
  \end{enumerate}

\end{lemma}

\begin{proof}
To go from (1) to (2), observe that
a nonabelian mixed Hodge structure always induces a mixed Hodge
structure on its Lie algebra compatible with bracket.  A Lie compatible mixed Hodge structure on
$\mathscr{G}_r(\pi)$ induces an Ind-MHS on its universal enveloping
algebra, compatible with the Hopf algebra structure.  This descends to $\Q[\pi]/J^{r+1}$ by (f) above.
A Hopf compatible mixed Hodge structure on  $\Q[\pi]/J^{r+1}$ induces one on $G_r(\pi)$ by
restriction.
\end{proof}

Hain \cite{hain} constructed a mixed Hodge structure on  the Hopf algebra $\Q[\pi_1(X,x)]/J^{r+1}$
which is equivalent (in the sense of the lemma) to the one constructed by Morgan on
$\mathscr{G}_r(\pi_1(X,x))$. These fit together to form an inverse
system as $r$ increases. In brief outline,
Chen had shown that  $\C\otimes \varprojlim \Q[\pi_1(X,x)]/J^{r+1}$ can be realized as the zeroth cohomology $H^0(B(x,\E^\dt(X),x))$ of a complex built from the $C^\infty$ de Rham complex
via the bar construction:
$$B^\dt(x,\E^*(X),x) =(\E^*(X))^{\otimes-\dt }$$
\begin{eqnarray*}
\pm d_B(\alpha_1\otimes \ldots \otimes \alpha_n) &=& i_x(\alpha_1)\alpha_2\otimes\ldots \otimes\alpha_n\\
&&+\sum (-1)^i\alpha_1\otimes\ldots \otimes \alpha_i\wedge\alpha_{i+1}\otimes\ldots \alpha_n\\
&&+ (-1)^ni_x(\alpha_n)\alpha_1\otimes\ldots \otimes\alpha_{n-1}\\
&& \pm d\alpha_1\otimes \alpha_2\otimes\ldots \otimes\alpha_n\\
&& + \ldots  
\end{eqnarray*}
where $i_x:\E^*(X)\to \C$ is the augmentation given by evaluation at $x$.
Hain showed how to extend $B(x,\E^\dt(X),x)$ to a cohomological mixed
Hodge complex (or more precisely a direct limit of such), and was thus able to  deduce the corresponding structure on cohomology. One thing that is more readily apparent in Hain's
approach is the dependence on base points. As $x$ varies, $\Q[\pi_1(X,x)]/J^{r+1}$
forms part of an admissible variation of mixed Hodge structure over $X$ called the tautological variation.
This is nontrivial since the monodromy representation is the natural conjugation homomorphism
$ \pi_1(X,x)\to Aut(\Q[\pi_1(X,x)]/J^{r+1})$.

Wojtkowiak \cite{wojt} gave a more algebro-geometric interpretation of Hain's construction, which
will be briefly described.
Bousfield and Kan defined the total space functor  $Tot$ \cite[chap X \S 3]{bk}, which is a kind of geometric realization, from the category of cosimplicial spaces to   the category of spaces. The image of the map
of cosimplicial schemes
$$
\xymatrix{
 X^{\Delta[1]}\ar@{}[r]^{=}\ar^{\pi}[d] & X\times X\ar@<1ex>[r]\ar@<-1ex>[r]\ar[d] & X\times X\times X\ar[l]\ar[d] & \ldots \\ 
 X^{\partial \Delta[1]}\ar@{}[r]^{=} & X\times X\ar@<1ex>[r]\ar@<-1ex>[r] & X\times X\ar[l] & \ldots
}
$$
under this functor is the  path space fibration $X^{[0,1]}\to X^{\{0,1\}}=X\times X$.
The horizontal maps on the top are diagonals (from left to right) or projections (from right to left); on the bottom they are all identities.
The total space of the  fibre $\pi^{-1}(x,x)$ is  the  space of loops of $X$ based at $x$, which
 is an $H$-space. Therefore $H^0(Tot (\pi^{-1}(x,x)),\Q)$ is naturally a Hopf algebra. This Hopf algebra,
which is described more precisely  in \cite{wojt}, can be identified with the coordinate ring of $\varprojlim G_r(\pi_1(X,x))$, or $G_r(\pi_1(X,x))$ if we truncate the cosimplicial space at the $r$th stage. This follows from  the  fact $H^0(Tot (\pi^{-1}(x,x)),\C)$ 
 can be computed using the total complex of the de Rham complex of the cosimplicial fibre, which is none other than $B(x,\E^\dt(X),x)$. Under this
identification the filtration by truncations induced on $H^0$ coincides
with the filtration by length of tensors on $B$.
 The  mixed Hodge structure on $H^0(Tot (\pi^{-1}(x,x)),\Q)$ can now 
be constructed using standard machinery:
take  compatible multiplicative mixed Hodge complexes
on each component of the cosimplicial space and then form the total complex
\cite[\S 5]{wojt}.  Furthermore
the tautological variations are given by
the $0$th total direct images of $\Q$ under $\pi|_X$ under the diagonal embedding $X\subset X\times X$.  A useful consequence of this point of view is that the MHS on
$\Q[\pi_1(X,x)]/J^{r+1}$ can be seen to come from a motive in Nori's sense \cite{cushman}.

Let $UMHS(X)$ ($U_rMHS(X)$) denote the subcategory of
unipotent admissible variations of mixed Hodge structure (with index
of unipotency at most $r$). The tautological variation associated to 
$\Q[\pi_1(X,x)]/J^{r+1}$ lies in $U_rMHS(X)$.
Given an object $V$ in $U_rMHS(X)$, the  monodromy representation
extends to an algebra homomorphism
 $$\Q[\pi_1(X,x)]/J^{r+1}\to End(V_x)$$
which is compatible with mixed Hodge structures.

\begin{thm}[Hain-Zucker {\cite{hz}}]\label{thm:hz}
  The above map gives an equivalence between $U_rMHS(X)$ and the category of 
  Hodge representations of $\Q[\pi_1(X,x)]/J^{r+1}$.
\end{thm}

The above equivalence respects the tensor structure.

  We note that every object
of $ U_rLoc(X)$ is a sum of subquotients of the local
system associated to the tautological representation 
$ \pi_1(X,x)\to Aut(\Q[\pi_1(X,x)]/J^{r+1})$. 
Therefore  $\phi(U_rMHS(X))= U_rLoc(X)$, where $\phi:U_rMHS(X)\to Loc(X)$
is the forgetful functor.
Consequently, we get a split exact sequence 
$$
\xymatrix{1\ar[r] &\pi_1(X,x)^{un_r}\ar[r] &  \Pi(U_rMHS(X))\ar[r] &
  MT\ar@<1ex>[l]\ar[r]& 1
}
$$ 
by theorem~\ref{thm:PiMT}.
In particular, $\pi_1(X,x)^{un_{r}}$ carries an NMHS which is a
quotient of the one on $\ph(X,x)$.

\begin{prop}\label{prop:morganhain}
  The above NMHS on $\pi_1(X,x)^{un_{r}}$ is equivalent to the Morgan-Hain structure on
  $\Q[\pi_1(X,x)]/J^{r+1}$.
\end{prop}

\begin{proof}
  By  lemma \ref{lemma:nmhspi1}, the Morgan-Hain structure on
  $\Q[\pi_1(X,x)]/J^{r+1}$ induces an NMHS on
  $\pi_1(X,x)^{un_{r}}$. Let $MH$ denote the semidirect product of
  $MT$ with this Hodge structure on $\pi_1(X,x)^{un_r}$. By theorem \ref{thm:hz}, we have a
  commutative diagram
 $$
  \xymatrix{ U_rLoc(X)\ar^{=}[d] &  U_rMHS(X)\ar[l]_{\phi}\ar^{\psi}[r]\ar^{\cong}[d] & MHS\ar^{\kappa}@<1ex>[l]\ar^{=}[d] \\
            Rep(\pi_1(X)^{un_r})&           Rep(MH)\ar[l]\ar[r]& Rep(MT)\ar@<1ex>[l]
} $$
where the functors $\psi$ and $\kappa$ are given by the fibre and the pullback along the constant map.
  Therefore $MH$ is   isomorphic to  $\Pi(U_rMHS(X))$ as a
  semidirect product.

  \end{proof}

  \begin{cor}
The Morgan-Hain NMHS on $\pi_1(X,x)^{un}$  is a quotient of $\ph(X,x)$.  
  \end{cor}

Hain has extended the above construction in  \cite{hain3}.
Given a representation $\rho:\pi_1(X,x)\to S$ to a reductive algebraic group,
the relative Mal\v{c}ev completion is the universal extension
$$1\to U\to \G\to S\to 1$$
of $S$ by a prounipotent group with a homomorphism $\pi_1(X,x)\to \G$ 
such that
$$
\xymatrix{\pi_1(X,x)\ar[r]\ar^{\rho}[dr]& \G\ar[d]\\
                     & S}
$$
commutes.
When $\rho:\pi_1(X,x)\to S=Aut(\V_x,\langle,\rangle)$ is the monodromy 
representation of a variation of Hodge structure with Zariski dense image,
Hain \cite{hain3} has shown that the relative Mal\v{c}ev completion carries a NMHS.

\begin{conj}\label{conj:relmal}
   The relative completion should carry an NMHS in general with $S$ equal
to the Zariski closure of $\pi_1(X,x)\to Aut(\V_x,\langle,\rangle)$. This should be a quotient of $\ph(X,x)$.
\end{conj}

I am quite confident about this. The essential point would construct an inner nonabelian 
variation of mixed Hodge structure on the family of $\G$ as the base point varies.
Then  proposition~\ref{prop:NVMHS} would give  a homomorphism 
$\ph(X,x)\to \G$. The main step  would be to establish an appropriate refinement of \cite[cor 13.11]{hain3}, and  I understand  that
 Hain, Matsumoto, Pearlstein, and Terasoma \cite{hmpt} have done this.
J. Pridham has pointed to me that his preprint \cite{pridham} may also have some
bearing on this conjecture.

\begin{rmk}\label{rmk:relmal}
For the applications given later in section~\ref{section:Kpi1}, only this weaker
statement on the existence of a morphism $\ph(X,x)\to \G$ extending
$\rho$ is needed.
  There is one notable case in which this can be deduced immediately.
If $\pi_1(X,x)$ is abelian, then $\G$ is necessarily abelian,
so it splits into a product $U\times S$. Morphisms $\ph(X,x)\to U$
and $\ph(X,x)\to S$ can be constructed directly from propositions \ref{prop:morganhain} and  \ref{prop:NVMHS}.
\end{rmk}

\section{Cohomology}

Fix a theory of enriched local systems $E$.
Let $G$ be a nonabelian $E$-structure.
The category of  representations of $\Pi(E(pt))\ltimes G$ 
is equivalent to the category of $E$-representations
of $G$. Given such a representation $V$, let $H^0(G,V) = V^G$
and $H^0(\Pi(E(pt))\ltimes G,V)= V^{\Pi(E(pt))\ltimes G}$ be the subspaces of invariants. 
The action of $\Pi(E(pt))$  on $V$ descends 
to an action on $H^0(G,V)$.  These
are left exact functors on the category of 
representations of $Rep_\infty(\Pi(E(pt))\ltimes G)$.
Since this category  has enough 
injectives (c.f. \cite[I 3.9]{jantzen}),
we can define higher derived functors $H^i(G, V)$ and $H^i(\Pi(E(pt)\ltimes G,V)$.
Note that $H^i(G, V)$ is an $E$-structure since it is  derived from a functor
from $Rep_\infty(\Pi(E(pt)\ltimes G)\to Rep_\infty(\Pi(E(pt))$.
Alternatively, $H^i(G, V)$  can  be computed as the cohomology of a bar or Hochschild complex $C^\dt(G, V)$  \cite[I 4.14]{jantzen}, which is a complex in
$ Rep_\infty(\Pi(E(pt))$.
We can define
cohomology of $\Pi(E(pt))$ by taking $G$ trivial. We observe that
these functors are covariant in $V$ and contravariant in $G$.

There are products 
$$H^i(G,V)\otimes H^j(G,V')\to H^{i+j}(G,V\otimes V')$$
compatible with $E$-structures. These can be constructed by either using
standard formulas for products on the complexes $C^\dt(G,-)$, or
by identifying $H^i(G, V)= Ext^i(\Q,V)$ and using the Yoneda pairing
$$
  Ext^i(\Q,V)\otimes Ext^j(\Q,V')\to  Ext^i(\Q,V)\otimes Ext^j(V,V\otimes V')
\to  Ext^{i+j}(\Q,V\otimes V')
$$
\begin{lemma}\label{lemma:HiG}
Given an $E$-representation  $V$ of $G$,
$H^i(G,V) = \varinjlim H^i(G/H_{j}, V^{H_j})$
where $H_j$ runs over all normal subgroups stable under $\Pi(E(pt))$ such
that $G/H_j$ is finite dimensional. 
\end{lemma}

\begin{proof}
 By lemma~\ref{lemma:NPHS} 
$G = \varprojlim G/H_j$ is an inverse
limit of algebraic groups with nonabelian $E$-structures.
Clearly $V = \varinjlim  V^{H_j}$ and so
$$H^0(G,V) = \varinjlim  H^0(G/H_j,V^{H_j})$$ 
The lemma follows from exactness of direct limits.
\end{proof}

\begin{lemma}
Fix a Hodge representation  $V$ of an NMHS $G$, and an
Ind-MHS $U$. Then
  \begin{enumerate}
\item $H^i(MT, U) =0$ for $i>1$. 
\item 
  We have an exact sequence
$$0\to H^1(MT, H^{i-1}(G,V))\to H^i(MT\ltimes G,V)\to  H^0(MT,
H^{i}(G,V))\to 0$$
  \end{enumerate}
\end{lemma}
\begin{proof}
Write $U$ as a direct limit of finite
dimensional Hodge structures $U_j$, then
 $$H^i(MT, U) = \varinjlim H^i(MT, U_j)= \varinjlim Ext_{MHS}^i(\Q, U_j) $$
Beilinson \cite{beilinson} shows that the higher $Ext$'s vanish, which
implies the first statement.

For the second statement, note that we can construct a
Hochschild-Serre spectral sequence 
$$E^{pq}_2 = H^p(MT, H^q(G, V))\Rightarrow H^{p+q}(MT\ltimes G,V)$$
in the usual way. This reduces to the given exact sequence thanks
to (1). 
\end{proof}

\begin{lemma}\label{lemma:HiGV}
Let $G$ be a nonabelian $E$-structure and
$V$ an $E$-representation.
  If $G$ is  pro-reductive then $H^i(G,V)=0$ for $i>0$.
In general 
$$H^i(G,V) = H^0(G^{red},H^i(U,V))$$
where $U$ is the pro-unipotent radical and $G^{red} = G/U$.
\end{lemma}

\begin{proof}
When $G$ is pro-reductive, its  category of representations  is
semisimple. Therefore $H^0(G,-)$ is exact. So higher cohomology
must vanish. 

In general, the Hochschild-Serre spectral sequence
$$E^{pq}_2 = H^p(G^{red}, H^q(U, V))\Rightarrow H^{p+q}( G,V)$$
will collapse to yield the above isomorphism.
\end{proof}

\begin{prop}\label{prop:Hipi}
Let $X$ be a smooth variety. Then for each object $V\in E(X)$:
\begin{enumerate}
\item $H^i(\pi^E_1(X,x),V)$ carries a canonical $E$-structure.
\item There is natural morphism $H^i(\pi^E_1(X,x),V)\to H^i(X,V)$ of
  $E$-structures.
\item  There are $\Q$-linear maps 
$$H^i(\pi^E_1(X,x), V)\to H^i(\pi_1(X,x), V)\to H^i(X, V)$$
whose composition is the map given in (2).
\end{enumerate}
 
\end{prop}

\begin{proof}
 $H^i(\pi_1^E(X,x),V)$ carries an $E$-structure by the discussion preceding lemma \ref{lemma:HiG}.
Note that $H^0(\pi_1^E(X,x), V) = V_x^{\pi_1^E(X)}$ is nothing but the
monodromy invariant part of $V_x$, and 
$H^i(\pi^E_1(X,x), V)$ is the universal $\delta$-functor extending
it, in the sense of \cite{groth}. By  axiom (E4), $H^i(X,V)$ with its $E$-structure 
also forms a $\delta$-functor, and there is an isomorphism
$$H^0(\pi_1^E(X,x),V)\cong H^0(X,V)$$
of $E$-structures. Therefore (2) is a consequence of universality.

After disregarding $E$-structure,
we can view $H^i(\pi^E_1(X,x), V)$ as a universal $\delta$-functor from
$E(X)\to Vect_\Q$.  Therefore, from the isomorphisms
$$H^0(\pi_1^E(X,x),V)\cong H^0(\pi_1(X,x),V)\cong H^0(X,V),$$
we deduce  $\Q$-linear maps
$$H^i(X,V)\leftarrow H^i(\pi_1^E(X,x),V)\to  H^i(\pi_1(X,x),V)$$
The leftmost map was the one constructed in the previous paragraph.
Since group cohomology is also a universal $\delta$-functor from the category
of $\Q[\pi_1(X,x)]$-modules to $Vect_\Q$,
we can complete this to a commutative triangle
$$
\xymatrix{
 H^i(\pi_1^E,V)\ar[r]\ar[d] & H^i(\pi,V)\ar@{-->}[ld] \\ 
 H^i(X,V) & 
}
$$

\end{proof}

\section{Hodge theoretic $K(\pi,1)$'s}\label{section:Kpi1}

The map 
\begin{equation}
  \label{eq:Hipi}
H^i(\pi_1^E(X,x), V)\to H^i(X, V)  \tag{*}
\end{equation}
constructed in the previous section is trivially an isomorphism
for $i=0$, but usually not  in general. For instance if $i=2$, $V=\Q$  and $X=\PP^1$, the cohomology 
groups are $0$ on the left
and $\Q$ on the right.

Let us say that $X$ is a $K(\pi^E,1)$ (or a Hodge theoretic
$K(\pi,1)$ when
$E=MHS$) if \eqref{eq:Hipi}  is an isomorphism for every $V\in E(X)$.
When $X$ is a $K(\pi,1)$ in the usual sense, then it is a  $K(\pi^E,1)$
if and only if 
\begin{equation}
  \label{eq:Hipi2}
H^i(\pi_1^E(X,x), V)\to H^i(\pi_1(X,x), V)  \tag{**}
\end{equation}
is an isomorphism for all $V\in E(X)$.
The following is a straight forward modification of  \cite[p 13, ex 1]{serre}.

\begin{lemma}\label{lemma:serre}
The following are equivalent
\begin{enumerate}
\item (\ref{eq:Hipi}) is an isomorphism for  all $i\le n$ and injective for $i=n+1$ for all $V\in $Ind-$E(X)$.
\item (\ref{eq:Hipi}) is an isomorphism for  all $i\le n$ and injective for $i=n+1$ for all $V\in E(X)$.
\item (\ref{eq:Hipi}) is surjective for all $i\le n$  and all $V\in E(X)$.
\item (\ref{eq:Hipi}) is surjective for all $i\le n$  and all $V\in \text{Ind-}E(X)$.
\item For all $V\in \text{Ind-}E(X)$, $1\le i\le n$, and $\alpha\in H^i(X, V)$, there exists $V'\in \text{Ind-}E(X)$ containing $V$ such that the image of $\alpha$ in $H^i(X, V') $ vanishes.
\end{enumerate}
In particular, $X$ is a  $K(\pi^E,1)$ if 
(\ref{eq:Hipi}) is surjective for all $i$ and $V\in E(X)$.
\end{lemma}

\begin{proof}
The implications (1)$\Rightarrow$(2)  and (2) $\Rightarrow$(3) are clear. For   (3)$\Rightarrow$(4), we can use
the fact that cohomology commutes with filtered direct limits.
 The implication (4)$\Rightarrow$(5) follows because Ind-$E(X)$ contains enough
injectives \cite{jantzen}. Any injective object $V'\supset V$ will satisfy the conditions of (5) assuming
(4).

Finally,  we prove  (5)$\Rightarrow$(1). This is the only nontrivial step.
An exact sequence
$$0\to V\to V'\to V'/V\to 0$$
yields a diagram
$$
\xymatrix{
 H^{i-1}(\pi_1^E,V')\ar[r]\ar^{f}[d] & H^{i-1}(\pi_1^E,V'/V)\ar[r]\ar^{g}[d] & H^i(\pi_1^E,V)\ar[r]\ar[d] & H^i(\pi_1^E,V')\ar[d] \\ 
 H^{i-1}(X,V')\ar[r] & H^{i-1}(X,V'/V)\ar[r] & H^i(X,V)\ar[r] & H^i(X,V')
}
$$

We first prove surjectivity of (\ref{eq:Hipi}) for $i\le n$ by induction.
This is trivially true when $i=0$, so we may assume $i>0$.
Given $\alpha\in H^i(X, V)$, we may choose $V'$ so that $\alpha$ has trivial image
in $H^i(X,V')$.
 Then $\alpha$ lifts to $ H^{i-1}(X,V'/V)$ and
hence to some $\beta\in H^{i-1}(\pi_1^E,V'/V)$. The image of
$\beta$ in $H^i(\pi_1^E,V)$ will map to $\alpha$ as required.

Now we prove injectivity of (\ref{eq:Hipi}) for $i\le n+1$ by induction.
We may assume that $i>0$. Let $V'$ be  injective in Ind-$E(X)$.
Suppose that 
$$\alpha\in \ker[H^i(\pi_1^E(X,x), V)\to H^i(X, V)] $$ 
Then $\alpha$ can be lifted to $\beta\in  H^{i-1}(\pi_1^E,V'/V)$.
Since the maps labelled $f$ and $g$ are isomorphisms, a simple
diagram chase shows that $\beta$ lies in the image of $H^{i-1}(\pi_1^E,V')$.
Therefore $\alpha=0$.
\end{proof}

\begin{prop}\label{prop:Kpi1}
When $E=MHS$, for all $V\in MHS(X)$ the map \eqref{eq:Hipi} 
is an isomorphism for  $i\le 1$ and injective for $i=2$
assuming conjecture~\ref{conj:relmal} holds. This is true unconditionally
if $\pi_1(X)$ is abelian.
\end{prop}

\begin{proof} 
As is well known, for any $V$ the map
$$H^i(\pi_1(X,x),V)\to H^i(X,V)$$
is an isomorphism for $i=1$ and injective for $i=2$.
Thus, by this remark and the previous lemma,
 it is enough to prove that the map \eqref{eq:Hipi2} to
group cohomology is surjective for  $i=1$.
By the 5-lemma and induction on the length of the weight filtration,
it is sufficient to the prove this when $V$ is pure.

Let $V$ be a variation of pure Hodge structure, and let 
$1\to U\to \G\to S\to 1$ be the 
associated relative Mal\v{c}ev completion. By lemma~\ref{lemma:HiGV},
$$H^1(\G,V) \cong H^0(S, H^1(U, V)) \cong Hom_S(U/[U,U], V)$$
By \cite[prop 10.3]{hain3}
$$H^1(\pi_1(X,x), V) \cong Hom_S(U/[U,U], V)$$
Therefore the natural map $H^1(\G,V)\to H^1(\pi_1,V)$ is an
isomorphism. As this factors through \eqref{eq:Hipi2}
(by conjecture~\ref{conj:relmal} or remark~\ref{rmk:relmal}),
\eqref{eq:Hipi2}  must be a surjection in degree $1$.
\end{proof}

\begin{thm}
A (not necessarily affine) connected commutative algebraic group is a Hodge theoretic $K(\pi,1)$.
Assuming conjecture~\ref{conj:relmal},
  a  smooth affine curve is a Hodge theoretic $K(\pi,1)$.
\end{thm}

\begin{proof}
Suppose that $X$ is a commutative algebraic group. Then it is homotopy equivalent to a torus.
In particular, it is a $K(\pi,1)$. So it suffices to check surjectivity of \eqref{eq:Hipi2}.
The group $\pi_1(X)$ is abelian and finitely generated,
which implies that the pro-algebraic completion  of $\pi_1(X)$ is a
commutative algebraic group. Therefore the same is true for $\ph(X)$.
After extending scalars, it follows that the group $\ph(X)\otimes \bar \Q$
is a product of $\mathbb{G}_a$'s, $\mathbb{G}_m$'s and a finite abelian group.
Consequently, any irreducible representation $V$ of $\ph(X)\otimes \bar \Q$ is
one dimensional. For such a module, the K\"unneth formula implies that
$$
H^i(\pi_1(X), V) = 
\begin{cases}
\wedge^i H^1(\pi_1(X), V) &\text{if $V$  is trivial}\\
0 & \text{otherwise}  
\end{cases}
$$
Therefore \eqref{eq:Hipi2} is surjective in this case.
By applying (an appropriate modification of) 
lemma~\ref{lemma:serre} to the category of semisimple $\ph(X)\otimes \bar \Q$ representations,
we can see  that  \eqref{eq:Hipi2} is an isomorphism.
We note that
$$
H^i(\ph(X,x), V\otimes \bar \Q)\cong H^i(\ph(X,x), V)\otimes \bar \Q
$$
$$
H^i(\pi_1(X,x), V\otimes \bar \Q)\cong H^i(\pi_1(X,x), V)\otimes \bar \Q
$$
and likewise for the map between them.
Thus we may extend scalars in order to test bijectivity in \eqref{eq:Hipi2}.
After doing so, we see that \eqref{eq:Hipi2} is an isomorphism by
an induction on the length of a Jordan-Holder series.

  Let $X$ be a smooth affine curve. This is a $K(\pi,1)$ with a free
  fundamental group. Since a free group has cohomological dimension one,
\eqref{eq:Hipi2}  is surjective in all degrees assuming \ref{conj:relmal},
by the previous proposition. 

\end{proof}

I expect that  Artin neighbourhoods are also Hodge theoretic $K(\pi,1)$'s.
This would give a large supply of such spaces.
Katzarkov, Pantev and Toen have established an analogous result in their setting \cite[rmk 4.17]{kpt2}.
Although, their proofs do not translate directly into the present framework, I
suspect that an appropriate modification may.



\begin{thebibliography}{9999}

\bibitem[A]{arapura} D. Arapura, {\em A category of motivic sheaves},
	ArXiv:0801.0261

\bibitem[B]{beilinson} A. Beilinson, {\em Notes on absolute Hodge
    cohomology}, Applications of algebraic K-theory to algebraic geometry and number theory, AMS (1987)

\bibitem[BK]{bk} A. Bousfield, D. Kan, {\em Homotopy limits, completions and localizations},
LNM 304, Springer-Verlag (1972)

\bibitem[C]{cushman} M. Cushman, {\em Morphisms of curves and the fundamental
group}, Contemp Trends in Alg. Geom. and Alg. Top., World Scientific (2002)

\bibitem[DM]{dm} P. Deligne, J. Milne, {\em Tannakian categories},
in LNM 900, Springer-Verlag (1982)



\bibitem[D1]{deligne1} P. Deligne, {\em Theorie de Hodge II},
 Inst. Hautes \'Etudes Sci. Publ. Math 40 (1971)




\bibitem[D2]{deligne2} P. Deligne, {\em Le groupes fondemental de la droit projective moin trois points}, Galois Groups over $\Q$, Springer-Verlag (1989)


\bibitem[D3]{deligne-tan} P. Deligne, {\em Cat\'egories Tannakiennes},
Grothendieck Festschrift, Birkhauser (1990)


\bibitem[G]{groth} A. Grothendieck, {\em Sur quelques points d alg\'ebre  homologiques}, 
Tohoku Math. J. 9 (1957)

\bibitem[H1]{hain} R. Hain, {\em The de Rham homotopy theory of complex algebraic varieties. I}, K-theory 3 (1987)





\bibitem[H2]{hain3} R. Hain, {\em The Hodge de Rham theory of relative Mal\v{c}ev
    completion,}  Ann. Sci. École Norm. Sup. 31 (1998), 47--92.

\bibitem[HZ]{hz} R. Hain, S. Zucker, {\em Unipotent variations of
mixed Hodge structure}, Inv. Math 88 (1987)


\bibitem[HMPT]{hmpt} R. Hain, M. Matsumoto, G. Pearlstein, T. Terasoma,
{\em Tannakian Fundamental Groups of Categories of Variations of Mixed
Hodge Structure}, in preperation

\bibitem[J]{jantzen} J. Jantzen, {\em Representations of
algebraic groups}, 2nd ed., AMS (2003)

\bibitem[K]{kashiwara} M Kashiwara, {\em A study of a variation of
mixed Hodge structure},  Publ. Res. Inst. Math. Sci. 22 (1986)

\bibitem[KPT1]{kpt1} L. Katzarkov, T. Pantev, B. Toen, {\em Schematic
homotopy types and non-abelian Hodge theory}, Compositio Math (2008)

\bibitem[KPT2]{kpt2} L. Katzarkov, T. Pantev, B. Toen, {\em 
Algebraic and topological aspects of the schematization functor},
arXiv.math.AG/0503418 


\bibitem[Ko]{kollar} J. Koll\'ar, {\em Shafarevich maps and automorphic forms},
Princeton U. Press (1995)

\bibitem[M]{morgan} J. Morgan, {\em The algebraic topology of
smooth algebraic varieties}, 
 Inst. Hautes \'Etudes Sci. Publ. Math 48 (1978)

\bibitem[P]{pridham} J. Pridham, {\em Formality and splitting of real non-abelian mixed Hodge 
structures,}  ArXiv:0902.0770


\bibitem[Q]{quillen} D. Quillen, {\em Rational homotopy theory},
  Ann. Math 90 (1969)

\bibitem[Sa1]{saito1} M. Saito, {\em Mixed Hodge modules and
    admissible variations}, CR Acad. Sci. Paris 309 (1989)

\bibitem[Sa2]{saito} M. Saito, {\em Mixed Hodge modules}, 
 Publ. Res. Inst. Math. Sci. 26 (1990)


\bibitem[Se]{serre} J.P. Serre, {\em Cohomologie Galois}, 5th ed., LNM 5,
 Springer-Verlag (1994)


\bibitem[Si]{simpson} C. Simpson, {\em Higgs bundles and local
systems,}  Inst. Hautes \'Etudes Sci. Publ. Math 75  (1992)



\bibitem[SZ]{sz} J. Steenbrink, S. Zucker, {\em Variation of
mixed Hodge structures I}, Inv. Math. 80 (1985)


\bibitem[W]{wojt} Z. Wojtkowiak, {\em Cosimplicial objects in algebraic
geometry}, Alg. K-theory and Alg. Top., Kluwer (1993)

\end{thebibliography}
\end{document}